\documentclass[reqno, 11pt]{amsart}  
\usepackage{amsmath}
\usepackage{amssymb}
\usepackage{hyperref}
\hypersetup{colorlinks=true}

\usepackage[left=1.125in,right=1.1in,top=1.1in,bottom=1.125in]{geometry}

\newtheorem{prop}{Proposition}[section]

\newtheorem{thm}[prop]{Theorem}

\numberwithin{equation}{section}

\newcommand{\shrink}[1]{ {\scriptstyle {\textstyle {#1} } } }
\newcommand{\smfrac}[2]{ \shrink{ \frac{#1}{#2} } }
\newcommand{\lin}{\langle} 
\newcommand{\rin}{\rangle} 
\newcommand{\nt}{\negthinspace}

\newcommand{\dist}{\mathrm{dist}} 
\newcommand{\Feas}{\mathrm{Feas} }

\newcommand{\xcomega}{X_{\Omega}}

\usepackage{array,multirow}
\usepackage{algpseudocode}
\usepackage{algorithm}

\begin{document}

\newpage 
$ \textrm{~} $ \quad \vspace{-3mm}

\title[Stochastic Convex Feasibility Problem]{A Different Perspective \\ On The Stochastic Convex Feasibility Problem}

\begin{abstract}
We investigate the stochastic convex feasibility problem, the problem of finding $ x $ satisfying a system of constraints $ f_{\omega}(x) \leq 0 $ ($ \omega \in \Omega $), where the functions $ f_{\omega} $ are convex, and where $ \Omega $ is the sample space in a probability triple $ (\Omega, {\mathcal F}, P) $. Our aim is to solve the problem approximately, in the sense that given values $ \epsilon > 0 $ and $ 0 < \Gamma < 1 $, a point $ x $ is found for which the probability measure $ P $ satisfies $ P \{ \omega \in \Omega \mid f_{\omega}(x) \leq   \epsilon \} \geq 1 -  \Gamma  $. Our primary vehicle for obtaining such a point $ x $ is a subgradient method taking ``Polyak steps'', for which the only input is the size of minibatches to be employed in the effort. We prove stochastic iteration bounds for the general setting, and improved bounds when a notion of ``H\"{o}lderian growth'' is satisfied, showing, for example, linear convergence when the growth is linear (as happens, say, for a finite system of linear inequalities). We also derive deterministic iteration bounds for knowing ``with high confidence'' that an iterate $ x_k $ is in hand which fulfills the desired inequality, $ P\{ \omega \in \Omega \mid f_{\omega}(x_k) \leq   \epsilon \}  \geq 1 -  \Gamma  $. 
\end{abstract} \vspace{-3mm}

\author[J. Renegar and S. Zhou]{James Renegar and Song Zhou}
\address{School of Operations Research and Information Engineering,
 Cornell University, U.S.} 
\thanks{Research supported NSF grant DMS-1812904.}


\vspace{-7mm}
  
\maketitle

\section{{\bf Introduction}} 

The convex feasibility problem -- the problem of finding a solution to a system of convex constraints -- is an integral component of convex optimization, and from some perspectives the two are equivalent, in the deterministic setting. In the stochastic setting the feasibility problem is less studied than the problem that often is referred to as ``stochastic optimization'', the problem of minimizing a function of the form $ \mathbb{E}_{\omega \in \Omega} [ f_{\omega}(x)  ] $, where $ \Omega $ is the sample space in a probability triple $ ( \Omega, {\mathcal F}, P ) $. 

For us the stochastic convex feasibility problem is that of finding, or approximating, a point $ x $ satisfying a system of functional constraints $ f_{\omega}(x) \leq 0 $ ($ \omega \in \Omega $) where the functions $ f_{\omega} $ are convex.\footnote{A close cousin is the problem of finding, or approximating, a point in the intersection of closed convex sets, $ \cap_{\omega \in \Omega} X_{\omega} $, for which (deterministic and randomized) projection methods have been extensively studied (see Necoara, Richt{\'a}rik and Patrascu \cite{necoara2019randomized} for discussion that fits well with our paper). The intersection problem can be (superficially) viewed as a special case of the functional constraint problem, by defining $ f_{\omega}(x) = \dist(x, X_{\omega}) $ (Euclidean distance from $ x $ to $ X_{\omega} $) -- then the full gradient step at $ x $ corresponds to projecting $ x $ onto $ X_{\omega} $. However, this viewpoint does not begin to capture the rich body of geometric ideas that have arisen and matured in the literature on projection methods, and so the intersection problem is best considered as a cousin rather than a special case of the functional constraint problem.} Beginning with Polyak \cite{polyak2001random} (and perhaps most recently, Necoara and Nedi\'c \cite{necoara2021minibatch}), the prominent approach for designing iterative methods is to view the feasibility problem, under mild assumptions, as being equivalent to the stochastic optimization problem in which the goal is to minimize $ \mathbb{E} [ f_{\omega}^+(x)  ] $ (where $ f_{\omega}^+(x) := \max \{ 0, f_{\omega_{\ell}}(x) \} $). Then if an algorithm employs minibatches in making iterations, it is natural to take a (sub)gradient step at the current iterate $ x_k $ for the sample-expectation function $ x \mapsto \frac{1}{L} \sum_{\ell = 1}^L f_{\omega_{\ell}}^+(x) $, where $ \omega_1, \ldots, \omega_{L} $ are independently drawn from $ (\Omega, {\mathcal F}, P) $.\footnote{The size of minibatches, $ L $, can be chosen ``large but not too large'', so as to fully yet efficiently utilize available parallelism in making the step. Ideally, the number of iterations is reduced by a factor of $ 1/L $.} A manner in which algorithmic complexity is measured is by the expected number of iterations sufficient to reach an iterate $ x_k $ for which $ \mathbb{E} [ f_{\omega}^+(x_k)  ]  \leq \epsilon $, where $ \epsilon > 0 $ is a parameter (perhaps chosen by the user).  

We take a different route. For us the minibatch problem is to take a subgradient step at the current iterate $ x_k $ for the sample-maximum function $ x \mapsto \max_{\omega_1, \ldots, \omega_L} f_{\omega_{\ell}}(x) $.  Our manner of measuring complexity is by the expected number of steps sufficient to reach an iterate $ x_k $ for which ``most'' of the functions satisfy $ f_{\omega}(x) \leq \epsilon $, in the sense that $ P \{ \omega \mid f_{\omega}(x_k) \leq \epsilon \} \geq 1 - \Gamma $, where $ 0 < \Gamma < 1 $ is a parameter signifying the proportion of functions that can be ignored. (The only input to the main algorithm is $ L $, the size of minibatches, not $ \Gamma $, $ \epsilon $, or any other parameter.)

This alternative approach to the stochastic convex feasibility problem results in surprisingly tractable analyses. For example, bounding the expected number of iterations for the main algorithm is done by computing the hitting time for a compound Bernoulli process, elementary probability. Besides bounding the expected number of iterations quite generally, we easily establish concentration bounds on the number of iterations, and more interesting, we establish much-improved bounds when a notion akin to ``H\"{o}lderian growth'' is satisfied, for all degrees of growth, not just the linear growth of piecewise-linear convex functions or the quadratic growth of strongly convex functions. Finally, we establish the analogous results under a slight modification to the algorithm which results in the user knowing with ``high confidence'' an iterate $ x_k $  is in hand that satisfies the desired inequality, $ P \{ \omega \mid f_{\omega}(x_k) \leq \epsilon \} \geq 1 - \Gamma $. Perhaps surprisingly, the iteration bounds here are deterministic -- all of the probability gets wrapped into the confidence level (albeit at the expense of potentially large minibatches). 

In the following three sections, a progression of three theorems is established. Immediately following the statement of a theorem, we provide discussion, including when appropriate, relations to the existing literature.

\section{{\bf Basics}} \label{basics} 

For a convex function $ h $ which attains its minimum value $ h^* $, the subgradient method with ``Polyak step size'' is the iterative algorithm given by $ x_+ = x - \frac{h(x) - h^*}{\| g \|^2} g $ where $ g \in \partial h(x) $ (subdifferential) and $ \| \, \, \| $ is the Euclidean norm. Under mild conditions, the iterates converge to the set of minimizers, $ \{ x \mid h(x) = h^* \} $. Likewise, for any constant $ c \geq h^* $, the iterative algorithm $ x_+ = x - \frac{\max \{ 0, h(x) - c \} }{\| g \|^2} g $ converges to  $ \{ x \mid h(x) \leq c \} $. We also refer to the steps made by this algorithm as being Polyak steps. 

Let $ X_{\Omega} := \{ x \mid f_{\omega}(x) \leq 0 \textrm{ for all $ \omega $}\} $, which we assume is nonempty. Input to the following algorithm is $ x_0 $, the initial iterate, and $ L $, the size of minibatches. Let $ z_0 $ denote the point in $ X_{\Omega} $ which is closest to $ x_0 $. We assume for all $ \omega $ that $ f_{\omega} $ is $ M $-Lipschitz\footnote{$ | f_{\omega}(x) - f_{ \omega}(y)| \leq M \| x - y \| $.} on an open neighborhood of the closed ball $ {\mathcal B} := \{ x \mid \| x - z_0 \| \leq \| x_0 - z_0 \| \} $, a consequence being that if $ x \in {\mathcal B} $ and $ g_{\omega} \in \partial f_{\omega}(x) $, then $ \| g_{\omega} \| \leq M $.  

	We are especially interested in the case where $ \Omega $ is finite but contains a vast number of elements, and the case that $ \Omega $ is infinite.  We assume samples are drawn independently from $ \Omega $ {\em  with replacement}. 
\vspace{2mm}

\noindent 
$ \textrm{~} $ \quad $ \mathtt{PolyakFM} $   (Polyak Feasibility Method) \\
$ \textrm{~} $ \quad 0) Input: $ x_0 $ (initial iterate), and integer $ L \geq 1 $ (size of minibatch) \\
$ \textrm{~} $ \quad  \qquad  Initialization: $ k = 1 $ \\
$ \textrm{~} $ \quad 1) Independently select samples $ \omega_{k,1}, \ldots, \omega_{k,L} $ from $ ( \Omega, {\mathcal F}, P) $. \\
$ \textrm{~} $ \quad 2) Determine $ \ell_k := \arg\max_{\ell = 1, \ldots, L} f_{k,\ell}(x_{k-1}) $, then let $ \epsilon_{k-1} := f_{k, \ell_k}(x_{k-1}) $ and compute
\[      x_k = \begin{cases} x_{k-1} & \textrm{if $ \epsilon_{k-1} \leq 0 $ } \\
                           x_{k-1} - \frac{ \epsilon_{k-1} }{\| g_{\omega_{k,\ell_k}} \|^2} \, g_{\omega_{k,\ell_k}} & \textrm{otherwise} \; ,    
\end{cases} \]
$ \textrm{~} $ \qquad  \qquad  \qquad \qquad  \qquad  \qquad  \qquad  \qquad   where  $ g_{\omega_{k,\ell_k}} \nt \nt \in \partial f_{\omega_{k,\ell_k}}(x_{k-1}) $. \\
$ \textrm{~} $ \quad 3) Let $ k+1 \rightarrow k $ and return to Step 1. 
\vspace{2mm}

To ease notation, let $ \bar{\omega}_k := (\omega_{k,1}, \ldots, \omega_{k,L}) $.  So as to view the algorithm as being entirely determined by $ x_0 $ and the sequence of sample vectors $ \bar{ \omega}_1, \bar{ \omega}_2, \ldots $, let us think of $ g_{\omega} = g_{\omega}(x) $ as being a fixed element of $ \partial f_{\omega}(x) $, that is, $ x \mapsto g_{\omega}(x) $ is a subgradient selection. This causes no loss of generality in that if $ g_{\omega_{k,\ell_k}}(x_{k-1}) $ is needed by the algorithm (i.e., if $ f_{\omega_{k,\ell}}(x_{k-1}) > 0 $ for some $ \ell $), the subgradient can be chosen on the fly, then used and immediately wiped from memory, because as will become clear, the point $ x_{k-1} $ will never again be visited by the algorithm, and hence $ g_{\omega_{k,\ell_k}}(x_{k-1}) $ will not be needed a second time. Likewise if index $ \ell_k $ is not uniquely determined, then choose any one of the candidate indices, use it, and discard from memory. Henceforth, we consider the algorithm to be entirely determined by $ x_0 $ and whatever sample vectors are drawn, $ \bar{\omega}_1, \bar{\omega}_2, \ldots $.   

For points $ x $ and scalars $ \epsilon > 0 $, let 
\[  \Omega (x,\epsilon) :=  \{ \omega \in \Omega \mid f_{\omega}(x) \leq  \epsilon \} \; . \]
  Given $ \epsilon $ and $ 0 < \Gamma < 1 $, the goal is to compute a point $ x $ for which $ P( \Omega (x, \epsilon)) \geq 1 -  \Gamma $. (Clearly, the algorithm itself is independent of $ \epsilon $ and $ \Gamma $.)
  
In the following theorem we assume $ 0 < \epsilon < M \dist(x_0, \xcomega) $ (Euclidean distance from $ x_0 $ to $ \xcomega $),  because if $ \epsilon \geq M \dist(x_0, \xcomega) $ then $ \Omega(x_0, \epsilon) = \Omega $, so $ x_0 $ achieves the goal of obtaining an iterate $ x_k $ satisfying $ P( \Omega(x_k,\epsilon) ) \geq 1 - \Gamma $.

\begin{thm} \label{thm.ba} 
For any scalars $ 0 < \Gamma < 1 $ and $ 0 < \epsilon < M \, \dist(x_0, \xcomega) $, the algorithm $ \mathtt{PolyakFM} $  reaches an iterate $ x_k $ satisfying $ P( \Omega(x_k,\epsilon)) \geq 1 - \Gamma $ within an expected number of iterations not exceeding 
\[  \mathbb{E}  =  \smfrac{1}{p} \left( \smfrac{M \, \dist(x_0, X_{\Omega})}{\epsilon } \right)^2 \quad \textrm{where } p = 1 - (1 - \Gamma)^L \; . \]
   Moreover,  we have the concentration inequalities 
\begin{equation}  \label{eqn.ba} 
k \geq  2 \mathbb{E}   \Rightarrow     \mathrm{Prob}\big( \, x_k \, \textrm{is first iterate satisfying } P(\Omega(x_k,\epsilon)) \geq 1 - \Gamma  \, \big) \, \leq \,  \smfrac{1}{2} \bigg(  \frac{1}{1 + \frac{1}{2} \frac{p}{1-p}}  \bigg)^{k - \lceil 2 \, \mathbb{E} \rceil}     \end{equation} 
\end{thm}

\noindent 
{\bf Remarks:} 
\begin{enumerate}

\item 
For integers $ L $ satisfying $ 1 \leq L \leq 1/ \Gamma $, it is not difficult to show $ \frac{1}{L \Gamma} \leq \frac{1}{p} <   \smfrac{2}{L \Gamma } $, and thus the theorem indicates that the expected number of iterations scales almost ideally with the size of minibatches.  

\item If $ \Omega $ is finite, say  $ \Omega = \{ 1, \ldots, m \} $, with $ P(i) = \frac{1}{m} $ for all $ i $, then the proof of the theorem can readily be modified to allow the minibatch samples to be drawn from $ \Omega $ without replacement (assuming $ 1 \leq L \leq m $). In this case the expected number of iterations scales perfectly, and choosing $ L = m $ results in a deterministic algorithm and iteration bound. However, if $ m $ is astronomical in size, choosing $ L \approx  m $ is computationally untenable. 

\item  In general it is impossible to know with certainty whether $ P( \Omega(x_k, \epsilon)) \geq 1 - \Gamma $. Nonetheless, ``high confidence'' embellishments to the algorithm can be developed, the subject of \S\ref{confidence}. There, the values $ \epsilon_{k} $ computed by $ \mathtt{PolyakFM} $ will play the fundamental role of being inferred values for which $ P( \Omega(x_k, \epsilon_k)) \geq 1 - \Gamma $ with high probability (but this depends on choosing $ L $ appropriately large, unlike the present section where $ L $ is allowed to be any positive integer). 

\item As for a number of papers on the feasibility problem, we can slightly generalize the setting by assuming the region of interest is $ X_{\Omega} \cap Y $, where $ Y $ is a closed convex set for which, given $ x $, it is easy to compute the orthogonal projection $ \Pi_{Y}(x) $, the closest point in $ Y $ to $ x $. The appropriate modification for $ \mathtt{PolyakFM} $ is simply to project onto $ Y $ the point computed in Step 2, that is, replace $ x_k $ by $ \Pi_Y(x_k) $.  As anyone familiar with analyses of subgradient methods can easily verify, all of our results remain valid for this slightly more general setting. Since our notation already is cluttered, we choose to leave $ Y $ out of the picture.

\end{enumerate}
 
The following simple proposition plays a significant role in the subgradient literature. The proposition allows us to reduce the proof of Theorem~\ref{thm.ba}  to establishing the expected time for a compound Bernoulli process to reach a specific value, that is, the proof reduces to establishing a hitting time. 

\begin{prop} \label{prop.bb} 
Assume $ h $ is a convex function, and assume $ x, z $ are points satisfying $ h(x) > 0 $ and  $ h(z) \leq 0 $. For $ g \in \partial h(x) $, the Polyak  step
\[  x_+ = x -  \smfrac{h(x)}{\| g \|^2} \, g  \;   \]
satisfies 
\begin{equation} \label{eqn.bb} 
  \| x_+ - z \|^2 \leq \| x - z \|^2 -  ( \smfrac{h(x)}{\| g \| } )^2 \; .  
  \end{equation}  
  \end{prop} 
\noindent {\bf Proof:} We have
\begin{align*} 
 \| x_+ - z \|^2 & = \| x - \smfrac{h(x)}{\| g \|^2} g - z \|^2 \\
                 & = \| x - z \|^2 + 2   \smfrac{h(x)}{\| g \|^2} \,  \lin g, z-x \rin + (  \smfrac{h(x)}{\| g \|} )^2 \\
   & \leq  \| x - z \|^2  + 2   \smfrac{h(x)}{\| g \|^2} (h(z) - h(x)) + (  \smfrac{h(x)}{\| g \|} )^2  \\
 & \qquad \textrm{(because $ g $ is a subgradient at $ x $)} \\
   & \leq  \| x - z \|^2 -   (  \smfrac{h(x)}{\| g \|} )^2 \\
&  \qquad \textrm{(using $ h(z) - h(x) \leq - h(x) $)} \; , 
\end{align*}    
 completing the proof.
 \hfill $ \Box $
\vspace{2mm}

\noindent {\bf Remark:} The proof is trivially extended to allow extrapolated subgradient steps, $ x_+ = x = \delta \smfrac{h(x)}{\| g \|^2} g $ for $ 0 < \delta < 2 $, the importance of which is that choosing $ \delta > 1 $ often results empirically in improved performance, and often results in obtaining a point in the interior of $ X_{\Omega} $ if the interior is not empty (we do not assume the interior is nonempty). For extrapolated steps, the inequality (\ref{eqn.bb})  becomes $  \| x_+ - z \|^2 \leq \| x - z \|^2 -  \delta (2 - \delta) ( \smfrac{h(x)}{\| g \| } )^2 $. This can be used exactly as we rely on (\ref{eqn.bb})  throughout the paper. We choose not to include this layer of detail so as to avoid the additional notation.
\vspace{2mm}

An easy corollary to the proposition is that for $ \mathtt{PolyakFM} $, we have $ \| x_{k} - z_0 \| < \| x_{k-1} - z_0 \| $ unless $ x_k = x_{k-1} $. In particular, all iterates lie in the ball $ {\mathcal B} $.  

Also an important consequence of the proposition is that the iterates of $ \mathtt{PolyakFM} $ satisfy
\begin{equation} \label{eqn.bc} 
 \dist(x_k, X_{\Omega})^2 \leq \dist(x_{k-1}, X_{\Omega})^2 - (\epsilon_{k-1}/M)^2 \quad \textrm{if $ \epsilon_{k-1} \geq 0 $}.\footnote{Indeed, in the proposition let $ h = f_{\omega_{k,\ell_k}} $, $ x = x_{k-1} $, and take $ z $ to be the point in $ X_{\Omega} $ that is closest to $ x $. Then (\ref{eqn.bb})  implies $ \| x_k - z \|^2 \leq \dist(x_k,X_{\Omega})^2 - (\epsilon_{k-1}/M)^2 $. Finally, use $ \dist(x_k, X_{\Omega}) \leq \| x_k - z \| $. } 
 \end{equation}

\noindent {\bf Proof of Theorem~\ref{thm.ba}:} We begin by inductively defining a variable which keeps track of whether a subgradient iterate $ x_k $ has been reached for which the desired inequality $   P( \Omega(x_k, \epsilon) ) \geq 1 - \Gamma $ holds -- in this case we rely on the symbol ``$ \bullet $'' -- and which also keeps track, when such an iterate has yet to be reached, of the number of iterations $ k $ for which $ \epsilon_{k-1} \geq  \epsilon  $. In this it is useful to define  $ j  + \bullet = \bullet $ for integers $ j $. 

Let 
\[   X_0 = \begin{cases}
 \bullet  & \textrm{if $   P( \Omega(x_0,\epsilon)) \geq 1 - \Gamma $} \\ 0 & \textrm{otherwise} \; , 
  \end{cases} \] 
and for sequences of vectors $ \bar{ \omega}_1, \bar{ \omega}_2, \ldots $ (where $ \bar{\omega}_k = (\omega_{k,1}, \ldots, \omega_{k,L}) $), inductively define
\[  X_k(\bar{\omega}_1, \ldots, \bar{ \omega}_k) = \begin{cases} \bullet & \textrm{if $ P( \Omega (x_k,\epsilon)) \geq 1 - \Gamma  $}  \\ X_{k-1} + 1 & \textrm{if $ P (\Omega (x_k,\epsilon)) < 1 - \Gamma $ and $ \epsilon_{k-1} \geq  \epsilon $} \\
X_{k-1} & \textrm{if $ P (\Omega (x_k,\epsilon)) < 1 - \Gamma $ and $ \epsilon_{k-1} <  \epsilon $}    \; , 
\end{cases}   \] 
where $ x_k = x_k(\bar{ \omega}_1, \ldots, \bar{ \omega}_k) $. (Note that if state $ \bullet $ is reached, it remains the state thereafter.)

The goal is to understand, probabilistically, when state $ \bullet $ is reached, assuming all samples $ \omega_{k,\ell}  $ are independently drawn from $ ( \Omega, {\mathcal F}, P) $.

According to Proposition~\ref{prop.bb}, the distances $ \| x_k - z_0 \| $ ($ k = 0, 1, \ldots $) form a decreasing sequence (strictly decreasing except when $ x_k = x_{k-1} $). Moreover, the number of indices $ k $ satisfying $ \epsilon_{k-1} \geq  \epsilon $ cannot exceed $ N :=  \lfloor ( \frac{ M \, \| x_0 - z_0 \|}{\epsilon } )^2 \rfloor $, else by induction, (\ref{eqn.bc})  implies $ \dist(x_k, \xcomega) $  becomes negative, a contradiction. Consequently, if indices $ k_1 < k_2 < \cdots < k_N $ all satisfy $ \epsilon_{k-1} \geq  \epsilon $, then necessarily $ \epsilon_{k_N +1}  <  \epsilon $ -- thus, necessarily $ \Omega(x_{k_N}, \epsilon) = \Omega   $ -- and hence $ X_{k_N}(\bar{ \omega}_1, \ldots, \bar{ \omega}_{k_N}) = \bullet $. Hence, any state occurring is either $ \bullet $ or one of the values $ 0, 1, \ldots, N-1 $.
 
Until reaching state $ \bullet $, the values $ X_k $ change in a simple stochastic manner expressed by the conditional probability
\begin{equation}  \label{eqn.bd}
     \mathrm{Prob}_{\omega_{k,1}, \ldots, \omega_{k,L}}( X_{k} = X_{k-1} + 1 \mid \bar{ \omega}_1, \ldots, \bar{ \omega}_{k-1} )  = 1 - P( \Omega (x_{k-1}, \epsilon ) )^L \geq 1 - (1 - \Gamma)^L = p \; . 
     \end{equation}  
Consequently, starting at $ x_0 $, the expected number of iterations sufficient to reach state $ \bullet $ is no greater than the expected number of steps for the compound Bernoulli process
\[ 
   Y_0 = 0, \quad Y_k = \begin{cases} Y_{k-1} + 1 & \textrm{with probability $ p $} \\
                                        Y_{k-1} & \textrm{with probability $ 1 - p $}   \end{cases} 
\]
to reach value $ N = \lfloor ( \frac{ M \, \| x_0 - z_0 \|}{\epsilon } )^2 \rfloor $, an expectation which is easily shown to equal $ \mathbb{E} :=   N/p  $. 

For the concentration inequalities, consider any integer $ k \geq 2 \mathbb{E}  $ and observe
\begin{align}
    \smfrac{  \mathrm{Prob}\big( Y_{k+1} = N \textrm{ and } Y_k = N-1 \big)}{ \mathrm{Prob}\big( Y_k = N \textrm{ and } Y_{k-1} = N-1 \big)} & = \smfrac{p \, {k \choose N-1} p^{N-1} (1-p)^{k-(N-1)}}{p \, {k - 1 \choose N-1} p^{N-1} (1-p)^{(k-1)-(N-1)}}  \label{eqn.be} \\
   & = \smfrac{k (1 - p)}{k - (N-1)  } \leq  \smfrac{\frac{2N}{p} (1-p)}{\frac{2N}{p} - (N-1)} 
    = \smfrac{1-p}{1-p + \frac{1}{2} p + \frac{p}{N}  } < \smfrac{1}{1 + \frac{1}{2} \frac{p}{1-p} } \; . \nonumber
\end{align}
The inequality (\ref{eqn.ba})  now follows by induction, the base case being $ k = \lceil 2 \mathbb{E} \rceil $ for which $ \mathrm{Prob}( Y_k = N \textrm{ and } Y_{k-1} = N-1 ) \leq \smfrac{1}{2}  $ (an immediate consequence of $ \mathbb{E} $ being an upper bound on the expected number of iterations).    \hfill $ \Box $

\section{{\bf Growth} } \label{growth}

In recent years a substantial literature has emerged that is focused on establishing improved iteration bounds of first-order methods when applied to to minimizing a convex function $ h $ possessing H\"{o}lderian growth, this being the property that there exist constants $ \mu > 0 $ and $ d \geq 1 $ such that for $ x $ in the domain of $ h $, we have the lower bound
\[    h(x) \geq \mu \, \dist(x, X^*)^d \]
where $ X^* $ is the set of minimizers of $ h $. When $ d = 1 $, this is sometimes referred to as ``linear'' growth or ``sharpness'', whereas for $ d =2 $, the function has ``quadratic'' growth. Strongly convex functions possess quadratic growth but so do many other convex functions. 

In this section we establish improved iteration bounds for $ \mathtt{PolyakFM} $   when the collection of functions $ f_{\omega} $ ($ \omega \in \Omega $) possesses a property exactly analogous to H\"{o}lderian growth.  For motivation, consider a collection of affine functions $ f_{\omega}(x) = \alpha_{\omega}^T x + \beta_{\omega} $ for which $ X_{\Omega}  $ is nonempty. Assuming $ \Omega $ is finite -- say, $ \Omega = \{ 1, \ldots, m \} $ -- a classic result of Hoffman \cite{hoffman2003approximate} asserts that there exists a constant $ \mu > 0 $ (now known as the Hoffman constant) such that for every $ x \notin X_{\Omega} $,  $ \max_{\omega} f_{\omega}(x) \geq \mu \, \dist(x, X_{\Omega} ) $, that is, there is at least one function $ f_{\omega} $ for which  $ f_{\omega}(x) \geq \mu \, \dist(x, X_{\Omega} ) $. 

However, when $ m $ is astronomical in size (and especially when $ \Omega $ is infinite), there is little use in the property of having, for each $ x $, at least one of the functions satisfying $ f_{\omega}(x) \geq \mu \, \dist(x,X_{\Omega}) $. Indeed, if there is only one such function $ f_{\omega} $ for an iterate $ x_k $, what are the chances that a randomized algorithm like $ \mathtt{PolyakFM} $ will make use of that function in computing $ x_{k+1} $? Instead, what is needed is that a ``considerable portion'' of samples $ \omega \in \Omega $ satisfy $ f_{\omega}(x) \geq \mu \, \dist(x,\Feas) $. This leads us to the following assumption. 
\begin{quote}
{\bf Growth Assumption:} There exist constants $ \mu > 0 $, $ d \geq 1 $ and $ 0 < \Delta  \leq  1 $  such that  for every $ x \in {\mathcal B} \setminus X_{\Omega} $, 
\[  P \{ \omega \in \Omega  \mid f_{\omega}(x) \geq \mu \, \mathrm{dist}(x, X_{\Omega})^d \} \geq  \Delta  \; . \]
\end{quote}
\vspace{3mm}

In the following theorem, as in Theorem~\ref{thm.ba}, we assume $ \epsilon < M \dist(x_0, \xcomega) $ since otherwise $ x_0 $ achieves the goal of obtaining an iterate $ x_k $ satisfying $ P( \Omega(x_k,\epsilon) ) \geq 1 - \Gamma $.

 \begin{thm}  \label{thm.ca} 
Make the Growth Assumption. For $ 0 < \epsilon < M \dist(x_0, \xcomega) $ and  $ 0 < \Gamma < \Delta $, $ \mathtt{PolyakFM} $   reaches an iterate $ x_k $ satisfying $ P( \Omega (x,\epsilon)) \geq 1 - \Gamma  $ -- and $ \dist(x_k,X_{\Omega}) \leq (\epsilon/\mu)^{1/d} $ -- within an expected number of iterations not exceeding 
\begin{equation}  \label{eqn.ca}
  \mathbb{E}' =   \smfrac{4}{p}   \left( 1 +  \left(  \smfrac{ M }{ \mu^{1/d} \epsilon^{1-1/d}} \right)^2 \min \left\{ \smfrac{1}{4^{1 - 1/d} - 1},  \log_2 \left( \smfrac{M \dist(x_0,X_{\Omega})}{\epsilon} \right)   \right\}   \right) \textrm{ where } p = 1 - (1 - \Gamma)^L \; . 
    \end{equation} 
 Moreover, (\ref{eqn.ba})  provides concentration inequalities if $ \mathbb{E} $ is replaced by $ \mathbb{E}' $, and the requirement $ k \geq 2 \mathbb{E} $  is replaced by $ k \geq 2 \mathbb{E}' $. 
\end{thm}

\noindent {\bf Remarks:} 
\begin{enumerate} 

\item When $ d = 1 $, the dependence on $ 1/\epsilon $ is only logarithmic. Necoara and Nedi\'c \cite{necoara2021minibatch} established a similar iteration bound for an algorithm which takes a subgradient step for the sample expectation, $ x \mapsto \frac{1}{L} \sum_{\ell = 1}^L f_{\omega_{\ell}}(x) $. Their growth assumption roughly translates as being that there exists $ \tilde{\mu} > 0  $ such that for all $ x \in {\mathcal B} $, there holds the lower bound $ \mathbb{E} [ f_{\omega}^+(x) ] \geq \tilde{\mu } \, \dist(x, X_{\Omega}) $.   While their minibatches require computing $ L $ subgradients, their bounds have no dependence on a quantity similar to $ 1/p $, and in fact, essentially replace $ 1/p $ by a value that under non-trivial conditions decreases proportionally to $ 1/L $, assuming that subgradient steps are extended to a particular length that in some important cases can be computed in practice. They do not consider growth of degree $ d > 1 $, nor do we see how their proofs might be extended to provide bounds of the quality of Theorem~\ref{thm.ca}, especially with regards to $ \epsilon $. Nonetheless, we would be remiss not to acknowledge that the interesting results of \cite{necoara2021minibatch} were the primary catalyst for our paper.

\item For all degrees of growth, the dependence on $ \epsilon $ improves on the worst-case complexity of subgradient methods in the deterministic setting, $ O(1/ \epsilon^2) $.  
 
\item With additional bookkeeping in the proof, the theorem can be improved if the growth assumption is satisfied both for, say, degree $ 1 $ with constant $ \mu_1 $, and for degree $ 2 $ with constant $ \mu_2 $. Then for $ x $ farther than distance $ \mu_1/\mu_2 $ away from $ X_{\Omega} $, the quadratic lower bound exceeds the linear lower bound, that is, the quadratic lower bound dominates. Assuming the quadratic growth dominates at $ x_0 $, the manner of proof of the theorem -- and the fact that the algorithm does not require $ d $ or $ \mu $ as input -- allow (\ref{eqn.ca})  to be replaced by a sum of two expressions, one obtained by substituting into (\ref{eqn.ca})  the values $ d = 2 $, $ \mu = \mu_2 $ and the choice $ \epsilon = \mu_1^2/\mu_2 $, and the other obtained by substituting $ d = 1 $, $ \mu = \mu_1 $ and replacing $ \dist(x_0, X_{\Omega}) $ by $ \mu_1/\mu_2 $. Similarly if the growth assumption is satisfied by any two degrees $ d < d' $ (or any three (or more) degrees, but then the notation becomes unwieldy). 

\end{enumerate}

\noindent {\bf Proof of Theorem~\ref{thm.ca}:} 
For non-negative integers $ i $, define $ e_i = (\frac{1}{2})^i M \dist(x_0,X_{\Omega}) $, and let $ I := \lceil \log_2 (\frac{M \dist(x_0,X_{\Omega})}{\epsilon}) \rceil $, the smallest integer $ i $ for which $ e_i \leq \epsilon $. 

Theorem~\ref{thm.ba}  implies for each $ i   $ that with probability 1, the algorithm reaches an iterate $ x_k $ for which $ P( \Omega (x_k, e_i) ) \geq 1 - \Gamma  $. Let $ x_{k_i} $ ($ = x_{k_i}(\bar{\omega}_1, \ldots, \bar{\omega}_{k_i}) $) be the first such iterate. Note that $ k_0 = 0 $ because $ \Omega(x_0, e_0) = \Omega $ due to the choice $ e_0 = M \, \dist(x_0, \xcomega) $. 

For $ i \geq 1 $, if in Theorem~\ref{thm.ba} we substitute $ x_{k_{i-1}} $ for $ x_0 $ and $ e_{i} $ for $ \epsilon $, we find that the expected number of iterations to reach $ x_{k_i} $ from $ x_{k_{i-1}} $ does not exceed $   \frac{1}{p} ( \frac{M \, \dist(x_{k_{i-1}}, X_{\Omega})}{e_{k_i} } )^2 $. Consequently, beginning at $ x_0 $, the expected number of iterations to reach an iterate $ x_k $ satisfying $ P( \Omega(x_k,\epsilon)) \geq 1 - \Gamma $ does not exceed       
\begin{equation} \label{eqn.cb} 
\smfrac{1}{p} \sum_{i=1}^{I} \left(   \smfrac{M \dist(x_{k_{i-1}}, X_{\Omega})}{e_{i} } \right)^2  \, \, = \, \,  \smfrac{1}{p} \left( 4 +  \sum_{i=2}^{I} \left(   \smfrac{M \dist(x_{k_{i-1}}, X_{\Omega})}{e_{i} } \right)^2 \right) .  
\end{equation}

Since for $ i \geq 1 $, the point $ x_{k_i} $ satisfies $ P( \Omega (x_{k_i}, e_{i})) \geq 1 - \Gamma  $, our assumption $ \Gamma < \Delta $ implies there exists $ \omega_i \in \Omega (x_{k_i}, e_{i}) $ for which $ f_{\omega_i}(x_{k_i}) \geq \mu \, \dist(x_{k_i}, \xcomega)^d $. Consequently, we have both
\begin{equation}  \label{eqn.cc}
  e_{i} \geq  f_{\omega_i}(x_{k_i}) \textrm{ and }  f_{\omega_i}(x_{k_i}) \geq  \mu \, \dist(x_{k_i}, X_{\Omega})^d, \textrm{ and thus }  \dist(x_{k_i}, X_{\Omega}) \leq ( e_{i}/\mu  )^{1/d}    \; . 
\end{equation} 
Since $ e_i = \smfrac{1}{2}  e_{i-1} $, this gives for $ i \geq 2 $ that
\[  \smfrac{ \dist(x_{k_{i-1}}, X_{\Omega})}{e_i} \leq   \smfrac{2}{e_{i-1}} \left( \smfrac{ e_{i - 1}}{\mu} \right)^{1/d}  =  \smfrac{2}{ \mu^{1/d} e_{i-1}^{1 - 1/d} } = \smfrac{2}{ \mu^{1/d} (2^{I-i} e_{I-1})^{1 - 1/d}} \leq   \smfrac{2}{ \mu^{1/d} (2^{I-i} \epsilon)^{1 - 1/d}}  \; . \]
Substituting into (\ref{eqn.cb})  shows the expected time to reach an iterate $ x_k $ satisfying $ P( \Omega(x_k,\epsilon)) \geq 1 - \Gamma  $ does not exceed
\[ 
   \smfrac{1}{p}   \left( 4 +  \left(  \smfrac{ 2 M }{ \mu^{1/d} \epsilon^{1-1/d}} \right)^2 \, \sum_{i=1}^{I-1} \left( \smfrac{1}{4^{1 - 1/d}} \right)^{I-i} \right)  \; . \]
Finally, observe 
\[ 
    \sum_{i=1}^{I-1} \left( \smfrac{1}{4^{1 - 1/d}} \right)^{I-i}  \leq  \min \{ I - 1, \sum_{j=1}^{\infty} \left(  \smfrac{1}{4^{1 - 1/d}} \right)^j \}  
    = \min \{ I - 1 , \smfrac{ 1 }{4^{1 - 1/d} - 1} \} \; , \] 
 completing the proof that $ \mathbb{E}' $ is an upper bound on the expected number of iterations to reach an iterate $ x_k $ satisfying $ P( \Omega(x_k, \epsilon)) \geq 1 - \Gamma   $. That the iterate also satisfies $ \dist(x_k,X_{\Omega}) \leq (\epsilon/\mu)^{1/d} $ is established exactly like (\ref{eqn.cc}).
 
 To verify the concentration inequalities, first observe we could have established that $ \mathbb{E}'  $ is an upper bound on the expected number of iterations by relying on first principles as in the proof of Theorem~\ref{thm.ba}, where the focus was on determining the expected number of iterations for the compound Bernoulli compound $ Y_0, Y_1, \ldots  $  to first reach value $ N = \lfloor (  \frac{M \| x_0 - z_0 \| }{ \epsilon } )^2 \rfloor $. In the present context the relevant target is instead 
 \begin{align*}
    N' & = \sum_{i=1}^I  \Bigl\lfloor ( \smfrac{M \dist(x_{k_{i-1}}, X_{\Omega})}{ e_i } )^2  \Bigr\rfloor \\
       & \leq  4 \left( 1 +  \left(  \smfrac{ M }{ \mu^{1/d} \epsilon^{1-1/d}} \right)^2 \min \left\{ \smfrac{1}{4^{1 - 1/d} - 1},  \log_2 \left( \smfrac{M \dist(x_0,X_{\Omega})}{\epsilon} \right)  \right\}   \right)  \\ & = p \,  \mathbb{E}' \; .   
\end{align*}
The desired concentration inequalities follow, by substituting $ N' $ for $ N $ in (\ref{eqn.be})  (and substituting $ \mathbb{E}' $ for $ {\mathbb E}  $).    
 \hfill $ \Box $
 \vspace{2mm}

For general convex functions possessing H\"{o}lderian growth, the trick of using summations like (\ref{eqn.cb})  in establishing improved iteration bounds was first  developed extensively in \cite{renegar2021simple}, albeit in a non-stochastic setting, yet in a manner that applies to a wide range of first order methods, not just subgradient methods.

 \section{{\bf Confidence}} \label{confidence} 
 
While $ \mathtt{PolyakFM} $   will, with probability 1, compute an iterate $ x_k $ satisfying $ P(\Omega(x_k, \epsilon )) \geq 1 - \Gamma $, in general it is impossible or impractical to be certain that such an iterate has been reached (an exception being when $ \Omega $ consists of only a modest number of samples). 
In this section we present a modification to $ \mathtt{PolyakFM} $  resulting in an algorithm which determines, with high confidence, that an iterate satisfying $ P( \Omega (x_k, \epsilon)) \geq 1 - \Gamma  $ is in hand. 

The only change to $ \mathtt{PolyakFM} $ is that now the user inputs $ x_0 $  and values $ 0 < \Gamma < 1 $, $ 0 < \alpha < 1 $, the latter being an allowable probability of error. (The algorithm computes the batch size based on $ \Gamma $ and $ \alpha $.) Here the iterate pairs $ (x_k, \epsilon_k) $ satisfy $ P(\Omega (x_k, \epsilon_k) \geq 1 - \Gamma $ ``with confidence $ 1 - \alpha $'', by which is meant that across all sequences of samples $ \bar{\omega}_1, \bar{\omega}_2, \ldots $, the probability that even one error is made does not exceed $ \alpha $ (an error occurs when a pair $ (x_k, \epsilon_k) $ fails to satisfy the inequality $ P( \Omega(x_k,\epsilon_k)) \geq 1 - \Gamma $). For $ \epsilon > 0 $, we establish bounds on the number of iterations needed to reach a pair $ (x_k, \epsilon_k) $ for which $ \epsilon_k \leq \epsilon $. Perhaps surprisingly, these iteration bounds are deterministic  -- all of the probability gets wrapped into the level of confidence, $ 1 - \alpha $ (but at the expense of large minibatches if $ \Gamma \approx 0 $).

Following is the modified algorithm, in which the batch size $ L $ slowly increases with the number of iterations. 
\vspace{1mm}

\noindent 
$ \textrm{~} $ \quad $ \mathtt{confidentPFM} $   (Confident Polyak Feasibility Method) \\
$ \textrm{~} $ \quad  0) Input: $ x_0 $, $ 0 < \Gamma < 1 $, $ 0 < \alpha < 1  $ (allowable probability of error).    \\
$ \textrm{~} $ \quad  \qquad  Initialization: $ k = 1 $ \\
$ \textrm{~} $ \quad 1) Compute $ L = \lceil  \frac{1}{\Gamma} \ln ( \frac{2 k^2 }{\alpha}) \rceil $ and independently select samples $ \omega_{k,1}, \ldots, \omega_{k,L} $ from $ ( \Omega, {\mathcal F}, P) $. \\
$ \textrm{~} $ \quad 2) Determine $ \ell_k := \arg\max_{\ell = 1, \ldots, L} f_{\omega_{k,\ell}}(x_{k-1}) $, 
 then let $ \epsilon_{k-1} := f_{\omega_{k,\ell_k}}(x_{k-1}) $ and compute
\[      x_k = \begin{cases} x_{k-1} & \textrm{if $ \epsilon_{k-1} \leq 0 $} \\
                           x_{k-1} - \frac{ \epsilon_{k-1} }{\| g_{\omega_{k,\ell_k}} \|^2} \, g_{\omega_{k,\ell_k}} & \textrm{otherwise} \; ,    
\end{cases} \]
$ \textrm{~} $ \qquad  \qquad  \qquad \qquad  \qquad  \qquad  \qquad  \qquad \qquad   where  $ g_{\omega_{k,\ell_k}} \nt \nt \in \partial f_{\omega_{k,\ell_k}}(x_{k-1}) $. \\
$ \textrm{~} $ \quad 3) Let $ k+1 \rightarrow k $ and return to Step 1. 
\vspace{2mm}

In the following theorem, as in Theorems~\ref{thm.ba}  and \ref{thm.ca}, we assume the desired accuracy $ \epsilon $ satisfies $ \epsilon < M \, \dist(x_0, \xcomega) $ since otherwise all samples $ \omega $ satisfy $ f_{\omega}(x_0) \leq \epsilon $ and thus $ \epsilon_0 $, computed in the first iteration, satisfies $ \epsilon_0 \leq \epsilon $.
\newpage

\begin{thm} \label{thm.da} 
The probability that $ {\tt confidentPFM} $ makes no errors is at least   $ 1 -  \alpha $.  

If $ 0 < \epsilon < M \dist(x_0, \xcomega) $, then $ {\tt confidentPFM} $  computes $ \epsilon_{k} $ satisfying $ \epsilon_k \leq \epsilon $ 
\begin{equation}  \label{eqn.da} 
 \textrm{within } \,  1 + \Bigl\lfloor \left( \smfrac{M \, \dist(x_0, X_{\Omega})}{\epsilon}\right)^2 \Bigr\rfloor \textrm{  iterations.}
 \end{equation} \label{eqn.db} 
Moreover, if  (1) the Growth Assumption holds, (2) $ \Gamma < \Delta $ and (3) $ {\tt confidentPFM} $ makes no errors, then the number of iterations does not exceed
\begin{equation} 
     5 + 4 \left(  \smfrac{ M }{ \mu^{1/d} \epsilon^{1-1/d}} \right)^2 \min \left\{ \smfrac{1}{4^{1 - 1/d} - 1},  \log_2 \left( \smfrac{M \dist(x_0,X_{\Omega})}{\epsilon} \right)   \right\}  \; . \end{equation} 
\end{thm}

\noindent 
{\bf Remarks:} 
\begin{enumerate}

\item We know of no similar general result in the literature on subgradient methods, although there certainly have been stronger high-confidence results in optimization settings focused on narrower classes of functions (e.g., the results for strongly convex stochastic composite optimization in the seminal work of Ghadimi and Lan \cite[\S4.2]{ghadimi2012optimal}).

\item The main shortcoming of the setup  is that the number of samples per iteration is proportional to $ 1/\Gamma $, a large number of samples if $ \Gamma \approx 0 $. On the other hand, due to the slow increase in sample size per iteration, the total number of samples drawn in iterations $ 1, \ldots,  k $ is $ O( \frac{k}{\Gamma} \log( \frac{k}{\alpha}) $, ideal in $ \alpha $ and off by at most a factor of $ \log(k) $ from being ideal in $ k $.  

\end{enumerate}
   
\noindent {\bf Proof of Theorem~\ref{thm.da}:} First we verify that the probability of error is no greater than $ \alpha $. Let $ x_{k-1} = x_{k-1}(\bar{\omega}_1, \ldots, \bar{\omega}_{k-1}) $ be an iterate (possibly $ k - 1 = 0 $). Let $ L_k = \lceil \frac{1}{\Gamma} \ln( \frac{2k^2}{\alpha}) \rceil $, the number of independent samples $ \omega_{k,\ell} $ drawn for iteration $ k $. Observe that by definition of $ \epsilon_{k-1} $, all of the samples satisfy $ \omega_{k,\ell} \in \Omega(x_{k-1},\epsilon_{k-1} ) $.

Assume the  algorithm makes an error on iterate $ k $, that is, assume $ P( \Omega(x_{k-1}, \epsilon_{k-1})) < 1 - \Gamma   $. Then the probability of having chosen all samples from $ \Omega(x_{k-1},\epsilon_{k-1} ) $ is
\[  P( \Omega(x_{k-1}, \epsilon))^{L_k} < (1 - \Gamma)^{L_k}  = (1 - \Gamma)^{ \lceil  \frac{1}{\Gamma} \ln (\frac{2 k^2}{\alpha}) \rceil } \leq \smfrac{\alpha }{2 k^2}   \]
-- the final inequality holds because an integer $ n $ satisfies $ (1 - \Gamma)^n \leq \frac{\alpha}{2k^2}  $ if and only if $   n \geq \smfrac{ \ln (\frac{\alpha}{2k^2} )}{ \ln ( 1 - \Gamma)} = \smfrac{\ln (2k^2/\alpha)}{ - \ln (1 - \Gamma)} $, and hence if $ n \geq \frac{1}{\Gamma} \ln ( \frac{2 k^2}{ \alpha}) $. Consequently, the probability of $ \mathtt{confidentPFM} $ ever making an error does not exceed $ \sum_{k=1}^{\infty} \frac{\alpha}{2k^2} < \alpha $. 
 
We next establish (\ref{eqn.da}).  This is easily accomplished. Indeed, whenever $ \epsilon_{k-1} > \epsilon $, we have from (\ref{eqn.bc})  that
\[    \dist(x_{k}, X_{\Omega})^2 < \dist(x_{k-1}, X_{\Omega})^2 - (\epsilon/M)^2 \; , \]
an amount of decrease that can occur for at most $ \lfloor ( \frac{M \, \dist(x_0, X_{\Omega})}{\epsilon})^2 \rfloor $ iterations. Since $ \epsilon_k $ is computed one iteration later than $ x_k $, (\ref{eqn.da})  is thus a valid upper bound.

Now assume the Growth Assumption holds and $ \Gamma < \Delta $. As in the proof of Theorem~\ref{thm.ca}, for integers $ i \geq 0 $ let $ e_i = (\frac{1}{2})^i M \dist(x_0, \xcomega) $. Define $ I = \lceil \frac{M \, \dist(x_0, \xcomega)}{\epsilon} \rceil $, the smallest integer $ i $ for which $ e_i \leq \epsilon $.   

For all $ i $ we know from (\ref{eqn.da})  that for some iterate $ x_k $, the value $ \epsilon_k $ computed in the following iteration satisfies $ \epsilon_k \leq e_i $. Let $ x_{k_i} $ be the first such iterate. Of course $ k_0 = 0 $ since $ \Omega(x_0,e_0) = \Omega $ due to $ e_0 = M \, \dist(x_0, \xcomega) $.  

For $ i \geq 1 $, if in (\ref{eqn.da})  we substitute $ x_{k_{i-1}} $ for $ x_0 $ and $ e_{i} $ for $ \epsilon $, we find that the number of iterations to reach $ x_{k_i} $ from $ x_{k_{i-1}} $ does not exceed $ \lfloor ( \frac{M \, \dist(x_{k_{i-1}}, \xcomega)}{e_i})^2 \rfloor $ (one iteration less than (\ref{eqn.da})  because $ x_{k_i} $ is computed one iteration before $ e_{k_i} $). Consequently, beginning at $ x_0 $, the number of iterations to reach an iterate $ x_k $ for which $ \epsilon_k \leq \epsilon $  does not exceed       
\begin{equation} \label{eqn.dc} 
 \sum_{i=1}^{I} \left(   \smfrac{M \dist(x_{k_{i-1}}, X_{\Omega})}{e_{i} } \right)^2  \, \, = \, \,  4 +  \sum_{i=2}^{I} \left(   \smfrac{M \dist(x_{k_{i-1}}, X_{\Omega})}{e_{i} } \right)^2  .  
\end{equation}
Observe that this deterministic iteration bound is the exact analogue of the expected iteration bound (\ref{eqn.cb})  -- only the factor $ 1/p $ is now missing. From this point onwards the proof follows that of Theorem~\ref{thm.ca}, {\em assuming $ \mathtt{confidentPFM} $ makes no errors}. In particular, for the iterates $ x_{k_i} $ we need that $ P( \Omega(x_{k_i}, \epsilon_{k_i})) \geq 1 - \Gamma $, implying $ P( \Omega(x_{k_i}, e_i)) \geq 1 - \Gamma $.

Thus, assuming $ \mathtt{confidentPFM} $ makes no errors, and proceeding exactly as in the proof of Theorem~\ref{thm.ca} beginning at (\ref{eqn.cb})   (but disregarding the factor $ 1/p $), we find that the number of iterations required to reach $ x_{k_I} $ from $ x_0 $ does not exceed
\[ 4   \left( 1 +  \left(  \smfrac{ M }{ \mu^{1/d} \epsilon^{1-1/d}} \right)^2 \min \left\{ \smfrac{1}{4^{1 - 1/d} - 1},  \log_2 \left( \smfrac{M \dist(x_0,X_{\Omega})}{\epsilon} \right)   \right\}   \right) \; . \]
Since $ \epsilon_{k_I} $ ($ \leq \epsilon $) is not computed until the following iteration, (\ref{eqn.db})  is thus a valid deterministic bound.  \hfill $ \Box $

\bibliographystyle{plain}
\bibliography{stochastic_convex_feasibility} 
  
\end{document}